\documentclass{crm-article}

\usepackage{epstopdf}
\usepackage{citehack}
\usepackage{algorithm}
\usepackage{algorithmic}
\usepackage{graphicx}
\usepackage{url}
\usepackage{hyperref}
\hypersetup{
  colorlinks=true,
  linkcolor=black,
  citecolor=black,
  urlcolor=blue
}
\floatname{algorithm}{Алгоритм}

\usepackage{xfrac}

\newcommand{\E}{\mathbb{E}}
\newcommand{\R}{\mathbb{R}}
\newcommand{\PP}{\mathbb{P}}

\newtheorem{definition}{Определение}
\newtheorem{theorem}{Теорема}
\newtheorem{lemma}{Лемма}
\newtheorem{corollary}{Следствие}

\begin{document}

\journalVol{10}

\journalNo{1}
\setcounter{page}{1}

\journalSection{МАТЕМАТИЧЕСКИЕ ОСНОВЫ И ЧИСЛЕННЫЕ МЕТОДЫ МОДЕЛИРОВАНИЯ}
\journalSectionEn{Mathematical modeling and numerical simulation}

\journalReceived{09.11.2020.}

\journalAccepted{xx.xx.2020.}

\UDC{519.85}
\title{Метод эллипсоидов для задач выпуклой стохастической оптимизации малой размерности\footnote{Статья была подготовлена в ходе проектной смены "Современные методы теории информации, оптимизации и управления"\;Сириус 2-23 августа 2020 г. и научно-образовательной школы-конференции "Управление. Информация. Оптимизация"\; Сириус 23-29 августа 2020 г.}}
\titleeng{Ellipsoid method for convex stochastic optimization in small dimension}

\thanks{Авторы выражают благодарность Гасникову Александру Владимировичу за идею для написания работы.}
\thankseng{Authors express gratitude to Alexander Gasnikov for the idea of the paper.}

\author[1,2]{\firstname{Е.\,Л.}~\surname{Гладин}}
\authorfull{Егор Леонидович Гладин}
\authoreng{\firstname{E.\,L.}~\surname{Gladin}}
\authorfulleng{Egor L. Gladin}
\email{gladin.el@phystech.edu}
\affiliation[1]{Национальный исследовательский университет «Московский физико-технический институт»,\protect\\ Россия, 141701, Московская облаcть,
г. Долгопрудный, 
Институтский пер., д. 9.}
\affiliationeng{National Research University Moscow Institute of Physics and Technology,\protect\\ 9 Institutskiy per., Dolgoprudny, Moscow Region, 141701, Russia}
\affiliation[2]{Сколковский институт науки и технологий,\protect\\ Россия, 121205, г. Москва, Большой бульвар, д. 30с1.}
\affiliationeng{Skolkovo Institute of Science and Technology,\protect\\ Bolshoy Boulevard 30, bld. 1, Moscow, 121205, Russia}
\author[1]{\firstname{К.\,Э.}~\surname{Зайнуллина}}
\authorfull{Карина Эдуардовна Зайнуллина}
\authoreng{\firstname{K.\,E.}~\surname{Zaynullina}}
\authorfulleng{Karina E. Zaynullina}
\email{zaynullina.ke@phystech.edu}



\begin{abstract}
В статье рассматривается задача минимизации математического ожидания выпуклой функции. Задачи такого вида повсеместны в машинном обучении, а также часто возникают в ряде других приложений. На практике для их решения обычно используются процедуры типа стохастического градиентного спуска. В нашей работе предлагается решать такие задачи с использованием метода эллипсоидов с минибатчингом. Алгоритм имеет линейную скорость сходимости, поэтому требует намного меньше итераций, чем стохастический градиентный спуск. Это подтверждается в наших экспериментах, исходный код которых находится в открытом доступе. Для получения линейной скорости сходимости метода не требуется ни гладкость, ни сильная выпуклость целевой функции. Таким образом, сложность алгоритма не зависит от обусловленности задачи. В работе доказывается, что метод эллипсоидов с наперёд заданной вероятностью находит решение с желаемой точностью при использовании минибатчей, размер которых пропорционален точности в степени -2.Это позволяет выполнять алгоритм параллельно на большом числе процессоров, тогда как возможности для батч параллелизации процедур типа стохастического градиентного спуска весьма ограничены. Несмотря на быструю сходимость, общее число вызовов оракула для метода эллипсоидов может получиться больше, чем для стохастического градиентного спуска, который неплохо сходится и при маленьком размере батча. Количество итераций метода эллипсоидов зависит от размерности задачи квадратично, поэтому метод подойдёт для относительно небольших размерностей.

\end{abstract}

\begin{abstracteng}
The article considers minimization of the expectation of convex function. Problems of this type often arise in machine learning and a number of other applications.
In practice, stochastic gradient descent (SGD) and similar procedures are often used to solve such problems. We propose to use the ellipsoid method with minibatching, which converges linearly and hence requires significantly less iterations than SGD. This is verified by our experiments, which are publicly available. The algorithm does not require neither smoothness nor strong convexity of target function to achieve linear convergence. We prove that the method arrives at approximate solution with given probability when using minibatches of size proportional to the desired precision to the power -2. This enables efficient parallel execution of the algorithm, whereas possibilities for batch parallelization of SGD are rather limited. Despite fast convergence, ellipsoid method can result in a greater total number of calls to oracle than SGD, which works decently with small batches.
Complexity is quadratic in dimension of the problem, hence the method is suitable for relatively small dimensionalities.

\end{abstracteng} 
\keywordeng{stochastic optimization, convex optimization, ellipsoid method, minibatching}

\maketitle

\paragraph{Введение}
В приложениях часто возникает задача оптимизации математического ожидания некоторой функции. Так, основной задачей машинного обучения является
подбор параметров модели
с помощью минимизации риска:
\begin{equation}\label{risk_min}
    \min_{\theta} \left\{ L(\theta) := \E_{(x,y)} l\left(f(\theta), x, y\right) \right\},
\end{equation}
где $l\left(f(\theta), x, y \right)$~--- функция потерь решающего правила $f(\theta)$ на паре объект-ответ $(x, y)$. Поскольку реальное распределение данных неизвестно, на практике часто осуществляется минимизация эмпирического риска:
\begin{equation}\label{ERM}
    \tilde{\theta}_{ERM}=\arg \min_{\theta} \frac{1}{m} \sum_{i=1}^{m} l\left(f(\theta), x_{i}, y_{i}\right),
\end{equation}
где $\left\{ (x_i, y_i) \right\}_{i=1}^n$~--- элементы выборки.
Вместо того, чтобы работать с готовой выборкой, можно решать задачу в онлайн-режиме, получая обучающие примеры последовательно.
Для этого можно использовать, например, стохастический градиентный спуск (SGD):
\begin{equation}\label{SGD}
    \theta \leftarrow \theta - \alpha \nabla_{\theta} l\left(f(\theta), \xi_{i}\right),\quad \text{где } \xi_{i} = (x_{i}, y_{i}),\ \alpha \text{ --- размер шага метода.}
\end{equation}

В стохастической оптимизации часто рассматриваются два подхода к определению точности приближённого решения $\tilde{\theta}$ задачи \eqref{risk_min} \cite{dvurechensky2018}. Первый из них связан с математическим ожиданием точности по функции. В этом случае $\tilde{\theta}$ называется $\varepsilon$-решением \eqref{risk_min} для $\varepsilon>0$, если $\E L(\tilde{\theta}) - L_* \leq \varepsilon$, где $L_*$~--- оптимальное значение в задаче \eqref{risk_min}, математическое ожидание берётся по случайности, возникающей в стохастическом алгоритме. Второй подход позволяет ограничить вероятность того, что значение функции в приближённом решении сильно отклоняется от оптимального. В этом случае $\tilde{\theta}$ называется $(\varepsilon, \beta)$-решением \eqref{risk_min} для $\varepsilon>0,\ \beta \in (0,1)$, если $\PP \left\{  L(\tilde{\theta}) - L_* > \varepsilon \right\} \leq \beta$.

Если функция $L(\theta)$ выпуклой, SGD позволяет получить $\varepsilon$-решение задачи \eqref{risk_min} за $N = O\left(\frac{1}{\varepsilon^2} \right)$ итераций. Если же мы хотим, чтобы $\tilde{\theta}_{ERM}$ из \eqref{ERM} было $\varepsilon$-решением задачи \eqref{risk_min}, то потребуется \cite{shapiro2014} размер выборки $m \sim \frac{n}{\varepsilon^2}$, где $n$~--- размерность $\theta$. Впрочем, этот результат можно также привести к $m \sim \frac{1}{\varepsilon^2}$, если регуляризовать целевую функцию в \eqref{ERM}, добавив к ней слагаемое $\sim \varepsilon \left\| \theta \right\|_2^2$. Тем не менее, найти точное решение \eqref{ERM}~--- часто сложная, а порой и невыполнимая задача.

Процедуры типа \eqref{SGD} можно осуществлять параллельно. В негладком случае поиск $(\varepsilon, \beta)$~-~решения \eqref{risk_min} удаётся производить лишь на $\Theta\left(\ln \left(\beta^{-1}\right)\right)$ машинах за счёт батч параллелизации \cite{dvurechensky2018}, что не даёт большого преимущества. Если же целевая функция является гладкой, методы типа SGD допускают батч параллелизацию на $\Theta\left( \frac{1}{\varepsilon^{3/2}} \right)$ процессорах \cite{woodworth2018}, то есть число (параллельных) итераций ускоренного метода стохастического градиентного спуска будет всего $O\left( \frac{1}{\sqrt{\varepsilon}} \right)$. В настоящей статье предлагается использовать метод эллипсоидов с минибатчингом для решения задачи \eqref{risk_min}, который имеет линейную скорость сходимости при параллелизации на $\widetilde{O} \left( \frac{1}{\varepsilon^{2}} \right)$ машинах. Метод подойдёт лишь для относительно небольших размерностей, поскольку число итераций зависит от размерности задачи квадратично.

\paragraph{Постановка задачи и полученные результаты}

Рассмотрим задачу
\begin{equation}\label{problem:min_f}
    \min_{x \in Q} \left\{ f(x) := \E_{\xi} f(x, \xi) \right\},
\end{equation}
где $Q \subseteq \mathbb{R}^n$~--- выпуклое компактное множество с непустой внутренностью, размерность $n$ относительно небольшая (до ста), $f(x)$~--- непрерывная выпуклая функция. Обозначим $\displaystyle D := \sup_{x, y \in Q} \left\| x - y \right\|,\ B := \sup_{x, y \in Q} \lvert f(x) - f(y) \rvert,\ \rho$~--- радиус некоторого евклидова шара, содержащегося в $Q$. Здесь и далее $\| \cdot \|$ означает евклидову норму. Будем считать, что стохастический оракул $\partial_x f (x, \xi)$ удовлетворяет для некоторого $\sigma > 0$ условию
\begin{equation}
    \mathbb{E_\xi} \exp \left( \frac{\left\|\partial_x f (x, \xi) - g \right\|^2}{\sigma^2} \right) \leq \exp (1),
\end{equation}
где $g$~--- некоторый субградиент $f$ в точке $x$ (обозначение: $g \in \partial f (x)$).
В настоящей статье предлагается использовать метод эллипсоидов с минибатчингом и доказывается, что он позволяет найти $(\varepsilon, \beta)$-решение задачи \eqref{problem:min_f} для $\varepsilon>0, \beta \in (0,1)$ после
\begin{equation*}
    N \leq 2 n^2 \ln \left( \frac{D B}{\rho \varepsilon} \right)
\end{equation*}
итераций при размере батча
\begin{equation*}
    r = \widetilde{O} \left( \frac{\sigma^2 D^2}{\varepsilon^2} \right),
\end{equation*}
где $\widetilde{O}(\cdot)$ означает $O(\cdot)$ с точностью до логарифмического по $\varepsilon^{-1}$ и $\beta^{-1}$ множителя.

\paragraph{Доказательство сходимости}
Нам потребуется следующее определение.
\begin{definition}\label{delta_subgradient}
Пусть $\delta \geq 0,\ Q \subseteq \R^n$~--- выпуклое множество, $f: Q \to \R$~--- выпуклая функция. Вектор $g \in \R^n$ называется $\delta$-субградиентом $f$ в точке $x \in Q$, если
\begin{equation*}
    f(y) \geq f(x) + \langle g, y-x \rangle - \delta\quad \forall y \in Q.
\end{equation*}
Множество $\delta$-субградиентов $f$ в точке $x$ обозначается $\partial_\delta f(x)$.
\end{definition}
Заметим, что $\delta$-субградиент при $\delta$=0 совпадает с обычным субградиентом. Приведённая ниже теорема устанавливает, что если в методе эллипсоидов (см., например, \cite{bubeck2015convex}) вместо субградиента использовать $\delta$-субградиент (алгоритм \ref{alg:ellipsoid}), то неточность $\delta$ не будет накапливаться от шага к шагу.
\begin{algorithm}[h!]
	\caption{Метод эллипсоидов с $\delta$-субградиентом для задачи $\min_{x \in Q} f(x)$}
	\label{alg:ellipsoid}
	\begin{algorithmic}[1]
		\REQUIRE Число итераций $N \geqslant 1$,  $\delta \geqslant 0$, шар $\mathcal{B}_{R} 	\supseteq Q$, его центр $c$ и радиус $R$.
		\STATE $\mathcal{E}_0 := \mathcal{B}_{R},\quad H_0 := R^2 I_n,\quad c_0 := c$.
		\FOR{$k=0,\, \dots, \, N-1$}
		    \IF {$c_k \in Q$}
		        \STATE $w_k := w \in  \partial_\delta f(x) $ 
		        \IF {$w_k = 0$}
		            \RETURN $c_k$ 
		        \ENDIF
		    \ELSE
		        \STATE $w_k := w$, где $w \neq 0$ таков, что $Q \subset \{ x \in \mathcal{E}_k: \langle w, x-c_k \rangle \leqslant 0 \}$
		    \ENDIF
		    \STATE $c_{k+1} := c_k - \frac{1}{n+1}\frac{H_k w_k}{\sqrt{w_k^T H_k w_k}}$ \\
		    $H_{k+1} := \frac{n^2}{n^2-1} \left( H_k - \frac{2}{n+1}\frac{H_k w_k w_k^T H_k}{w_k^T H_k w_k} \right)$ \\
		    $\mathcal{E}_{k+1} := \{x: (x-c_{k+1})^T H_{k+1}^{-1} (x-c_{k+1}) \leqslant 1 \}$
		\ENDFOR
		\ENSURE $x^N = \arg\min\limits_{x \in \{c_0, ..., c_N \} \cap Q } f(x)$
	\end{algorithmic}
\end{algorithm}
\begin{theorem}{\cite{gladin2020spp}}\label{th:ellipsoidz}
    Пусть $Q$~--- компактное выпуклое множество, которое содержится в некотором евклидовом шаре радиуса $R$ и включает некоторый евклидов шар радиуса $\rho,\ $ $f: Q \rightarrow \mathbb{R}$~--- непрерывная выпуклая функция, число $B>0$ таково, что $|f(x) - f(x')| \leq B\ \forall x, x' \in Q$.
    После $N \geq 2n^2 \ln \frac{R}{\rho}$ итераций метод эллипсоидов с $\delta$-субградиентом (алгоритм \ref{alg:ellipsoid}) возвращает такую точку $x^N \in Q$, что
    $$
        f(x^N) - \min_{x \in Q} f(x_*) \leq \frac{B R}{\rho} \exp \left(-\frac{N}{2n^2} \right)+\delta.
    $$
\end{theorem}
Далее приведены две леммы, которые позволяют связать $\delta$-субградиент со стохастическим оракулом, использующим минибатчинг, для которого мы используем обозначение
\begin{equation*} 
    \overset{r}{\partial_x} f \left( {x, \{ \xi^l \}_{l=1}^r} \right) := \frac{1}{r} \sum_{i=1}^r \partial_x f (x,  \xi^l ),\quad \xi^l \text{~--- независимые реализации случайной величины } \xi.
\end{equation*}
\begin{lemma}\label{subgrad_norm}
    Пусть $g \in \partial f (x)$, и пусть для вектора $\tilde{g}$ выполнено $\left\| \tilde{g} - g \right\| \leq \varepsilon$, тогда $\tilde{g} \in \partial_{\delta} f(x)$ для $\delta = \varepsilon D$, где $D:= \sup_{x, y \in Q} \left\| x - y \right\|$.
\end{lemma}
\textit{Доказательство}\quad По неравенству Коши-Буняковского, для любого $y \in Q$ справедливо
\begin{equation}\label{cauchy}
    \langle \tilde{g} - g, y - x \rangle \leq \left\| \tilde{g} - g \right\| \cdot \left\| y - x \right\| \leq \left\| \tilde{g} - g \right\| \cdot D \leq \varepsilon D
\end{equation}
В силу выпуклости $f$ для $g \in \partial f (x)$ имеем
\begin{equation}\label{subgrad}
    f(y) \geq f(x) + \langle g, y - x \rangle\quad \forall y \in Q.
\end{equation}
Сложив \eqref{cauchy} и \eqref{subgrad}, получим
\begin{equation}\label{delta_subgrad}
    f(y) \geq f(x) + \langle \tilde{g}, y - x \rangle - \varepsilon D\quad \forall y \in Q,
\end{equation}
т.е. $\tilde{g} \in \partial_{\delta} f(x),\ \delta = \varepsilon D$.

\begin{lemma}\label{large_div_lemm}
    Пусть стохастический оракул $\partial_x f (x, \xi)$ удовлетворяет условию $\mathbb{E_\xi} \exp \left( \frac{\left\|\partial_x f (x, \xi) - g \right\|^2}{\sigma^2} \right) \leq~\exp \left(1\right)$, где $g \in \partial f (x)$,  тогда для любого $\beta \in (0, 1)$
    
    \begin{equation*}
        \mathbb{P} \left\{ \left\|\overset{r}{\partial_x} f \left( {x, \{ \xi^l \}_{l=1}^r} \right) - g \right\| < \left[\sqrt{2} + \sqrt{6\ln{\beta^{-1}}}\right] \cdot \frac{\sigma}{\sqrt{r}} \right\} \geq 1 - \beta \quad \forall x \in Q.
    \end{equation*}
\end{lemma}
\textit{Доказательство}\quad Согласно теореме 2.1(ii) из \cite{juditsky2008large}, для любого $\gamma \geq 0$,
\begin{equation*}
\mathbb{P} \left\{ \left\| S_r \right\| \geq \left[\sqrt{2} + \sqrt{2}\gamma\right] \cdot \sqrt{\sum_{i=1}^r \sigma_i^2} \right\} \leq \exp \left({-\frac{\gamma^2}{3}} \right),
\end{equation*} где  $S_r$~--- сумма независимых случайных векторов $\{ \zeta_i \}_{i=1}^r$ с нулевым математическим ожиданием, для которых выполнено $\mathbb{E} \exp \left( \frac{\left\| \zeta_i \right\|^2}{\sigma_i^2} \right) \leq \exp \left(1\right)$.
В нашем случае $\zeta_i = \partial_x f (x,  \xi^i ) - g$  и $ \left\|S_r\right\| =  r \cdot~\left\|\overset{r}{\partial_x} f \left( {x, \{ \xi^l \}_{l=1}^r} \right) - g \right\|$, поэтому 
\begin{equation}\label{inequality_batching}
\mathbb{P} \left\{ \left\|\overset{r}{\partial_x} f \left( {x, \{ \xi^l \}_{l=1}^r} \right) - g \right\| \geq \left[\sqrt{2} + \sqrt{2}\gamma\right] \cdot \sqrt{\frac{\sigma^2}{r}} \right\} \leq \exp \left({-\frac{\gamma^2}{3}} \right) 
\end{equation} 

Обозначим $\beta = \exp\left({-\frac{\gamma^2}{3}}\right)$, тогда $\gamma = \sqrt{3 \ln\beta^{-1}}$ и неравенство (\ref{inequality_batching}) приобретает вид:

\begin{equation*}
\mathbb{P} \left\{ \left\|\overset{r}{\partial_x} f \left( {x, \{ \xi^l \}_{l=1}^r} \right) - g \right\| \geq \left[\sqrt{2} + \sqrt{6 \ln\beta^{-1}}\right] \cdot \frac{\sigma}{\sqrt{r}} \right\} \leq \beta. 
\end{equation*} 

Тогда 
\begin{multline*}
\mathbb{P} \left\{ \left\|\overset{r}{\partial_x} f \left( {x, \{ \xi^l \}_{l=1}^r} \right) - g \right\| < \left[\sqrt{2} + \sqrt{6 \ln\beta^{-1}}\right] \cdot \frac{\sigma}{\sqrt{r}} \right\} = 
        \\ 
        = 1 - \mathbb{P} \left\{ \left\|\overset{r}{\partial_x} f \left( {x, \{ \xi^l \}_{l=1}^r} \right) - g \right\| \geq \left[\sqrt{2} + \sqrt{6 \ln\beta^{-1}}\right] \cdot \frac{\sigma}{\sqrt{r}} \right\}
        \geq 1 -\beta
\end{multline*}

    
\begin{corollary}\label{cons_batching}
Пусть стохастический оракул $\overset{r}{\partial_x} f \left( {x, \{ \xi^l \}_{l=1}^r} \right)$ был вызван $N$ раз, $\tilde{g}_i$~--- его значение на шаге $i \in \overline{1, N}$ в точке $x^i$. Для любого $\beta \in (0, 1)$ справедливо
\begin{equation*}
    \mathbb{P} \left( \bigcap\limits^N_{i=1} \left\{ \tilde{g}_i \in \partial_{\delta} f (x^i) \right\} \right) \geq 1 - \beta N, \quad \text{где } \delta = \left[\sqrt{2} + \sqrt{6\ln{\beta^{-1}}}\right] \cdot \frac{\sigma D}{\sqrt{r}}.
\end{equation*}
\end{corollary}
\textit{Доказательство}\quad Согласно леммам \ref{subgrad_norm} и \ref{large_div_lemm}, для любых $x \in Q,\, \beta \in (0, 1)$
\begin{equation*}
    \mathbb{P} \left\{ \overset{r}{\partial_x} f \left( {x, \{ \xi^l \}_{l=1}^r} \right) \in \partial_{\delta} f (x) \right\} \geq 1 - \beta, \quad \text{где } \delta = \left[\sqrt{2} + \sqrt{6\ln{\beta^{-1}}}\right] \cdot \frac{\sigma D}{\sqrt{r}}.
\end{equation*}
Это эквивалентно
\begin{equation*}
    \mathbb{P} \left\{ \overset{r}{\partial_x} f \left( {x, \{ \xi^l \}_{l=1}^r} \right) \notin \partial_{\delta} f (x) \right\} \leq \beta.
\end{equation*}
Значит, вероятность того, что хотя бы на одном из $N$ шагов стохастический оракул вернёт вектор, который не будет являться $\delta$-субградиентом, не превышает $\beta N$. Следовательно, вероятность того, что стохастический оракул вернёт $\delta$-субградиент на всех $N$ шагах, составляет не менее $1 - \beta N$.

\begin{theorem}
    Метод эллипсоидов со стохастическим градиентом возвращает $(\varepsilon, \beta)$-решение задачи \eqref{problem:min_f} для $\varepsilon>0, \beta \in (0,1)$ после
    \begin{equation*}
        N \leq 2 n^2 \ln \left( \frac{D B}{\rho \varepsilon} \right)
    \end{equation*}
    итераций при размере батча
    \begin{equation*}
        r = \widetilde{O} \left( \frac{\sigma^2 D^2}{\varepsilon^2} \right).
    \end{equation*}
\end{theorem}
\textit{Доказательство}\quad Согласно теореме \ref{th:ellipsoidz},  $\varepsilon$-решение задачи \eqref{problem:min_f} может быть получено за
\begin{equation*}
    N \leq 2 n^2 \ln \left( \frac{D B}{\rho \varepsilon} \right)
\end{equation*}
итераций метода эллипсоидов с $\frac{\varepsilon}{2}$-субградиентом (алгоритм \ref{alg:ellipsoid}). Воспользуемся следствием \ref{cons_batching} для определения размера батча $r$, необходимого для того, чтобы стохастический оракул возвращал $\frac{\varepsilon}{2}$-субградиент на каждой из $N$ итераций с вероятностью не менее $1 - \beta$:
\begin{equation*}
    \frac{\varepsilon}{2} = \left[\sqrt{2} + \sqrt{6\ln \frac{N}{\beta}}\right] \cdot \frac{\sigma D}{\sqrt{r}} \Longrightarrow r = \widetilde{O} \left( \frac{\sigma^2 D^2}{\varepsilon^2} \right).
\end{equation*}

\paragraph{Численный эксперимент}
Рассмотрим модель логистической регрессии для задачи классификации:
\begin{equation*}
    \hat{p}_x (w) = \frac{1}{1 + e^{-\langle w, x \rangle}},
\end{equation*}
где $x$~--- вектор признаков для объекта обучающей выборки (включая константный признак), $w$~--- веса модели. В качестве функции потерь выступает кросс-энтропия:
\begin{equation*}
    f_x (w) = y \ln \hat{p}_x (w) + (1-y) \ln \left(1 - \hat{p}_x (w) \right),
\end{equation*}
где $y \in \{0,1 \}$~--- класс объекта обучающей выборки. Задача оптимизации имеет вид:
\begin{equation*}
    \min_w \left\{ f(w) := \E_x f_x (w) \right\}
\end{equation*}

\begin{figure}[H]
    \centering
\includegraphics[scale=0.8]{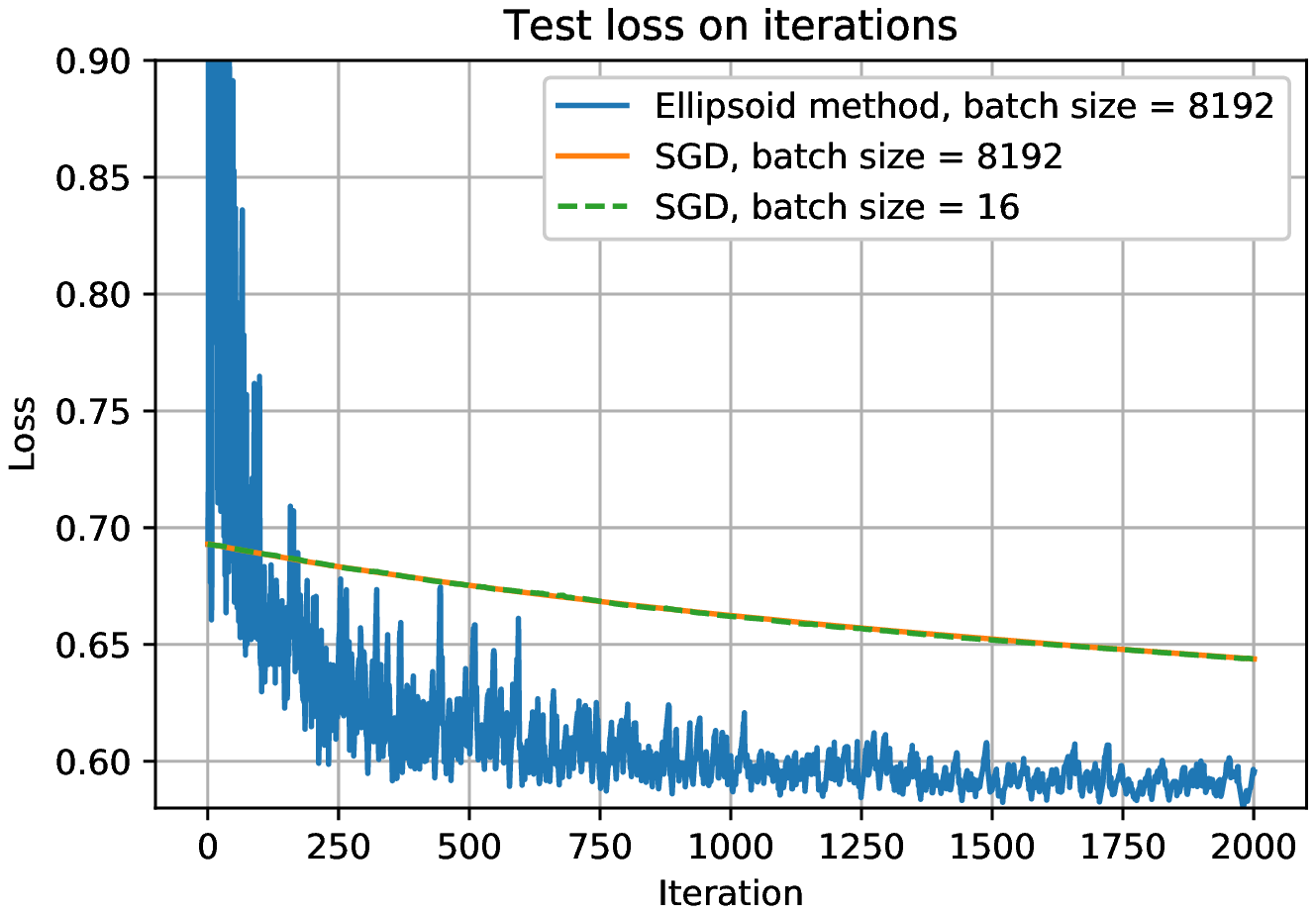}
    \caption{Сравнение метода эллипсоидов и стохастического градиентного спуска.}
    \label{paint1}
\end{figure}

В ходе экспериментов было проведено сравнение метода эллипсоидов и стохастического градиентного спуска при размере датасета $\sim 500 000$ объектов и размерности пространства признаков 55. Зависимость функции потерь на тестовой выборке от количества итераций для модели, обученной с помощью метода эллипсоидов с минибатчингом и стохастического градиентного спуска с минибатчингом отражена на графике \ref{paint1}. 

Как видно из графика, метод эллипсоидов сходится существенно быстрее стохастического градиентного спуска. 
Заметим, что метод эллипсоидов требует большого размера батча и будет эффективен только в случае, когда есть возможность осуществлять вычисления параллельно. Несмотря на быструю сходимость, общее число вызовов оракула для метода эллипсоидов может получиться больше, чем для SGD, который неплохо сходится и при маленьком размере батча.

Реализацию метода эллипсоидов на PyTorch и эксперименты можно найти в \href{https://github.com/Karina1997/math-optimization-ellipsoids-method}{GitHub репозитории}.

\begin{thebibliography}{99}

	
	
	
\bibitem[Гладин и др., 2020]{gladin2020spp}
	\textit{Гладин~Е.\,Л.,
    Курузов~И.\,А., Стонякин~Ф.\,С., Пасечнюк~Д.\,А., Алкуса~М.\,С., Гасников~А.\,В.} Решение сильно выпукло-вогнутых композитных седловых задач с небольшой размерностью одной из групп переменных. // arXiv preprint arXiv:2010.02280~--- 2020.
	
	\vspace{0.1cm}{\footnotesize{\it Gladin~E.} Reshenie sil'no vypuklo-vognutyh kompozitnyh sedlovyh zadach s nebol'shoj razmernost'yu odnoj iz grupp peremennyh. [Solving strongly convex-concave composite saddle point problems with a small dimension of one of the variables] // arXiv preprint arXiv:2010.02280~--- 2020.\par}
	
\bibitem[Bubeck, 2018]{bubeck2015convex}
	\textit{Bubeck S.} Convex optimization: algorithms and complexity. // Foundations and Trends in Machine Learning. – 2015. – V. 8. – №. 3--4. P. 231--357.
	
\bibitem[Dvurechensky et al., 2018]{dvurechensky2018}
	\textit{Dvurechensky~P.~E., Gasnikov~A.~V., Lagunovskaya~A.~A.} Parallel algorithms and probability of large deviation for stochastic convex optimization problems. // Numerical Analysis and Applications. – 2018. – V. 11. – №. 1.~–~P.33-37.

\bibitem[Juditsky et al., 2008]{juditsky2008large}
 	\textit{Juditsky~A., Nemirovski~A.} Large deviations of vector-valued martingales in 2-smooth normed spaces  //  arXiv preprint arXiv:0809.0813~--- 2008.
	
\bibitem[Shapiro et al., 2014]{shapiro2014}
	\textit{Shapiro~A., Dentcheva~D., Ruszczynski~A.} Lecture on stochastic programming. Modeling and theory. //  MPS-SIAM series on Optimization~--- 2014.
	
\bibitem[Woodworth et al., 2018]{woodworth2018}
	\textit{Woodworth~B. et al.} Graph Oracle Models, Lower Bounds, and Gaps for Parallel Stochastic Optimization // Advances in Neural Information Processing Systems. – 2018. – C. 8496-8506
	

	

	
	
	

	

	





	
	
	



	

	
\end{thebibliography}
\end{document}